\documentclass[11pt]{article}
\usepackage{latexsym,amssymb}
\usepackage{amsmath}

\newtheorem{Theorem}{\bf Theorem}[section]
\newtheorem{Lemma}{\bf Lemma}[section]
\newtheorem{Proposition}{\bf Proposition}[section]
\newtheorem{Corollary}{\bf Corollary}[section]
\newtheorem{Remark}{\bf Remark}[section]
\newtheorem{Example}{\bf Example}[section]
\newtheorem{Definition}{\bf Definition}[section]

\newenvironment{theorem}{\begin{Theorem}$\!\!\!$}{\end{Theorem}}
\newenvironment{lemma}{\begin{Lemma}$\!\!\!$}{\end{Lemma}}
\newenvironment{proposition}{\begin{Proposition}$\!\!\!$}{\end{Proposition}}
\newenvironment{corollary}{\begin{Corollary}$\!\!\!$}{\end{Corollary}}
\newenvironment{remark}{\begin{Remark}$\!\!\!$}{\end{Remark}}

\newenvironment{definition}{\begin{Definition}$\!\!\!$}{\end{Definition}}

\numberwithin{equation}{section}

\begin{document}

\title{Parabolic power concavity\\ 
and  parabolic boundary value problems}
\author{Kazuhiro Ishige\\
Mathematical Institute, Tohoku University\\
Aoba, Sendai 980-8578, Japan
\\
\\
and
\\
\quad
\\
Paolo Salani\\
Dipartimento di Matematica ``U. Dini'', Universit\`a di Firenze\\
viale Morgagni 67/A, 50134 Firenze}
\date{}
\maketitle

\pagestyle{myheadings}
\markright{Parabolic power concavity}
\begin{abstract}
This paper is concerned with power concavity properties of the solution to the parabolic boundary value problem
\begin{equation}
\tag{$P$}
\left\{
\begin{array}{ll}
\partial_t u=\Delta u +f(x,t,u,\nabla u) & \mbox{in}\quad\Omega\times(0,\infty),\vspace{3pt}\\
u(x,t)=0 & \mbox{on}\quad\partial \Omega\times(0,\infty),\vspace{3pt}\\
u(x,0)=0 & \mbox{in}\quad\Omega,
\end{array}
\right.
\end{equation}
where $\Omega$ is a bounded convex domain in ${\bf R}^n$ 
and $f$ is a nonnegative continuous function in $\Omega\times(0,\infty)\times{\bf R}\times{\bf R}^n$. 
We give a sufficient condition 
for the solution of $(P)$ to be parabolically power concave in $\overline{\Omega}\times[0,\infty)$. 
\end{abstract}
\vspace{10pt}
\section{Introduction}
We are concerned with the boundary value problem 
\begin{eqnarray}
\label{eq:1.1}
 & & \partial_t u=\Delta u +f(x,t,u,\nabla u)\qquad\mbox{in}\quad D,\vspace{3pt}\\
\label{eq:1.2}
 & & u(x,t)=0\qquad\qquad\qquad\qquad\,\,\, \mbox{on}\quad\partial D,
\end{eqnarray}
where $\partial_t=\partial/\partial t$, $D:=\Omega\times(0,\infty)$, 
$\Omega$ is a bounded convex domain in ${\bf R}^n$ $(n\ge 1)$, 
and $f$ is a nonnegative continuous function in $D\times{\bf R}\times{\bf R}^n$. 
In this paper we study power concavity properties  of the solution of problem~\eqref{eq:1.1} and \eqref{eq:1.2} 
with respect to the space and the time variables.  
For instance, it is shown that 
for any $\alpha\ge 1/2$, 
the solution $u$ of 
\begin{equation}
\label{eq:1.3}
\partial_t u=\Delta u +1\quad\mbox{in}\quad D,
\qquad
u(x,t)=0\quad\mbox{on}\quad\partial D
\end{equation}
is $\alpha$\,-parabolically $(1/2)$\,-concave in $\overline{D}$, 
that is, the function $v(x,t):=\sqrt{u(x,t^{1/\alpha})}$ is concave 
with respect to the variables $(x,t)\in\overline{D}$. 
\vspace{5pt}

Let us recall the notion of $p$\,-concavity for nonnegative functions, 
where $-\infty\le p\le\infty$. 
For $a$, $b>0$, $\lambda \in (0,1)$, and $p\in [-\infty, \infty]$,
we define
$$
M_p(a,b;\lambda)=
 \begin{cases}
\left[(1-\lambda) a^p+\lambda b^p\right]^{1/p}
& \mbox{if}\quad p\not\in\{-\infty, 0, \infty\},\\
a^{1-\lambda}b^{\lambda} & \mbox{if}\quad p=0,\\
\max\{a,b\} & \mbox{if}\quad p=\infty,\\
\min\{a,b\}& \mbox{if}\quad p=-\infty,
\end{cases}
$$
which is the $p$\,-(weighted) mean of $a$ and $b$ with ratio $\lambda$. 
Furthermore, for $a$, $b\ge 0$, we define $M_p(a,b;\lambda)$ as above if $p\ge 0$ and
$$
M_p(a,b;\lambda)=0\quad\mbox{if}\quad p<0\quad\text{and}\quad a\cdot b=0.
$$
\begin{definition}
\label{Definition:1.1}
Let $K$ be a convex set in ${\bf R}^m$ and $p\in[-\infty,\infty]$.
A nonnegative function $v$ defined in $K$ is said {\em $p$\,-concave} if
$$
v((1-\lambda)x+\lambda y)\geq M_p(v(x),v(y);\lambda)
$$
for all $x$, $y\in K$ and $\lambda\in(0,1)$. 
In the cases $p=0$ and $p=-\infty$, 
$v$ is also said log-concave and quasi-concave in $K$, respectively. 
\end{definition}
Notice that $p=1$ corresponds to usual concavity.

It follows from the Jensen inequality that 
if $v$ is $p$\,-concave in a convex set $K$, 
then $v$ is $q$\,-concave in $K$ for any $q\le p$ (see also \eqref{eq:2.2}). 
This means that quasi-concavity is the weakest concavity property one can imagine. 
\vspace{5pt}

Concavity of solutions for elliptic and parabolic partial differential equations is a classical subject 
and has fascinated many mathematicians. 
Especially in the elliptic case the literature is very large
and we just refer to the classical monograph by Kawohl \cite{Kaw} 
and the papers \cite{ALL}, \cite{BG}--\cite{BS}, \cite{cuoghis}, \cite{J}, \cite{Kenni}, \cite{Korev},  
some of which are closely related to this paper and the others include recent developments in this subject. 
Among (and before) others, 
Kennington~\cite{Kenni} improved the convexity maximum principle by Korevaar~\cite{Korev}
and established power concavity theorems for nonnegative solutions to elliptic boundary value problems
\begin{equation}
\label{eq:1.4}
\left\{
\begin{array}{ll}
\Delta u +f(x,u,\nabla u)=0 & \mbox{in}\quad\Omega\times(0,\infty),\vspace{3pt}\\
u=0 & \mbox{on}\quad\partial \Omega\times(0,\infty),
\end{array}
\right.
\end{equation}
under a suitable structure condition on the inhomogeneous term $f$. 
For readers's convenience, we recall the following result from \cite{Kenni} (see also \cite{ALL, BS, J} for some generalizations). 
\begin{theorem} {\rm \cite[Theorem~4.2]{Kenni}}. 
\label{Theorem:1.1} 
Let $\Omega$ be a bounded convex domain in ${\bf R}^n$ 
and $f$ a nonnegative function in $\Omega$. 
Let $u\in C^2(\Omega)\cap C(\overline{\Omega})$ satisfy 
\begin{equation}
\label{eq:1.5}
\left\{
\begin{array}{ll}
\Delta u +f(x)=0 & \mbox{in}\quad\Omega\times(0,\infty),\vspace{3pt}\\
u=0 & \mbox{on}\quad\partial \Omega\times(0,\infty). 
\end{array}
\right.
\end{equation}
\begin{itemize}
  \item[{\rm (i)}] 
  If $f$ is $q$\,-concave in $\Omega$ for some $q\ge 1$, 
  then $u$ is $p$\,-concave in $\overline{\Omega}$ with $p=q/(1+2q)$;
  \item[{\rm (ii)}] 
  If $f$ is a positive constant in $\Omega$, then $\sqrt{u}$ is concave in $\overline{\Omega}$. 
\end{itemize}
\end{theorem}
Compared with the elliptic equations, 
much less is known about the concavity properties for parabolic equations, 
and most of the results concern concavity properties with respect to the spatial variable only 
(see e.g. \cite{AI, BL, GK, IS1, IS2, KL}, \cite{Lee}--\cite{LV2} and references therein). 
Due to the concavity properties of the heat kernel and the first Dirichlet eigenfunction for Laplacian, 
log-concavity seems to be a natural property for the heat flow. 
Indeed, it is known that 
not only the log-concavity of the initial datum is preserved by the heat flow (see \cite{BL}, \cite{GK}, and \cite{Korev}), 
but also the solution of the heat equation in ${\bf R}^n$ becomes spatially log-concave in finite time 
provided that the initial function is nonnegative and has a compact support (see \cite{LV}). 
Notice that weaker concavity properties than log-concavity are not necessarily preserved by the heat flow 
(see \cite{IS1} and \cite{IS2}). 
On the other hand, 
inspired by \cite{Borell}, 
the authors of this paper 
introduced in \cite{IS3} and \cite{IS4}  the notions of parabolic and $\alpha$\,-parabolic quasi-concavity, 
and studied quasi-concavity properties involving the space and the time variables jointly 
for particular parabolic boundary value problems in a convex ring. 
See also \cite{HM} and \cite{PorruSerra} for results related to
space-time convexity of solutions of parabolic problems.
\vspace{5pt}

In this paper, following \cite{IS3} and \cite{IS4}, 
we introduce the notion of $\alpha$\,-parabolic $q$\,-concavity for nonnegative functions in a convex cylinder
and study parabolic power concavity properties of solutions for parabolic boundary value problems. 
\begin{definition}
\label{Definition:1.2}
Let $K$ be a convex set in ${\bf R}^n$, $Q:=K\times(0,\infty)$ and 
$\alpha,p\in[-\infty,\infty]$. 
A nonnegative function $v$ defined in $Q$ is said $\alpha$\,-parabolically $p$\,-concave if 
\begin{equation}
\label{eq:1.6}
v\big((1-\lambda)x_1+\lambda x_2, M_\alpha(t_1,t_2;\lambda)\big)\ge M_p\big(v(x_1,t_1),v(x_2,t_2);\lambda\big)
\end{equation}
for all $(x_1,t_1)$, $(x_2,t_2)\in Q$ and $\lambda\in(0,1)$. 
In particular, if $v$ is $(1/2)$\,-parabolically $p$\,-concave in $Q$, then 
it is simply said parabolically $p$\,-concave in $Q$. 
\end{definition}
Similarly to Definition~\ref{Definition:1.1},  
$\alpha$\,-parabolic log-concavity and $\alpha$\,-parabolic quasi-concavity correspond 
to $\alpha$\,-parabolic $0$\,-concavity and $\alpha$\,-parabolic $(-\infty)$\,-concavity, respectively. 
Roughly speaking, for some $\alpha\in{\bf R}\setminus\{0\}$, $v$ is $\alpha$\,-parabolically $p$\,-concave in $Q$ if 
\begin{itemize}
  \item $v$ is a constant function in $Q$ for $p=\infty$;
  \item $v(x,t^{1/\alpha})^p$ is concave in $Q$ for $p>0$;
  \item $\log v(x,t^{1/\alpha})$ is concave in $Q$ for $p=0$;
  \item $v(x,t^{1/\alpha})^p$ is convex in $Q$ for $p<0$;
  \item the superlevel sets $\{(x,t)\in Q:v(x,t^{1/\alpha})>\mu\}$ are convex for every $\mu\ge 0$ for $p=-\infty$.
\end{itemize}
Obviously, if a function $v$ is $\alpha$\,-parabolically $p$\,-concave for some $\alpha\in[-\infty,+\infty]$, 
then $v(\cdot,t)$ is spatially $p$\,-concave at any fixed time $t$.
\vspace{5pt}

Now we are ready to state a result on parabolic power concavity
for the boundary value problem~\eqref{eq:1.1} and \eqref{eq:1.2}. 
Here we focus on the case where $f$ depends only on the space variable.    
Theorem~\ref{Theorem:1.2} is an extension of Theorem~\ref{Theorem:1.1} to parabolic equations 
and a typical application of the main theorem of this paper, 
which is stated in Section~3. 
\begin{theorem}
\label{Theorem:1.2}
Let $\Omega$ be a bounded convex domain in ${\bf R}^n$, $D:=\Omega\times(0,\infty)$, 
and $f$ a nonnegative function in $\Omega$. 
Let $u\in C^{2,1}(D)\cap C(\overline{D})$ satisfy 
\begin{equation}
\label{eq:1.7}
\left\{
\begin{array}{ll}
\partial_tu=\Delta u +f(x) & \mbox{in}\quad D,\vspace{3pt}\\
u=0 & \mbox{on}\quad\partial D. 
\end{array}
\right.
\end{equation}
\begin{itemize}
  \item[{\rm (i)}]  
  If $f$ is $q$\,-concave in $\Omega$ for some $q\ge 1$, 
  then $u$ is $\alpha$\,-parabolically $p$\,-concave in $\overline{D}$ with $p=q/(1+2q)$ and $\alpha\ge 1/2$;
 \item[{\rm (ii)}] 
 If $f$ is a positive constant in $\Omega$, then $\sqrt{u(x,t^{1/\alpha})}$ is concave in $\overline{D}$ for any $\alpha\ge 1/2$.
\end{itemize}
\end{theorem}
As an application of Theorem~\ref{Theorem:1.2}, 
apart from the corresponding spatial power concavity of $u(\cdot, t)$ at any fixed time $t$,
we obtain the following power concavity properties in time of the heat energy 
$$
H(t):=\int_\Omega u(x,t)dx 
$$
associated to the source $f$ in $\Omega$. 
\begin{corollary}
\label{Corollary:1.1}
Assume the same conditions as in Theorem~{\rm\ref{Theorem:1.2}}. 
\begin{itemize}
  \item[{\rm (i)}] 
  If $f$ is $q$\,-concave in $\Omega$ for some $q\ge 1$, 
  then $H(t)$ is $\beta$\,-concave in $(0,\infty)$ with $\beta=q/[(n+2)q+1]$;
  \item[{\rm (ii)}]  
  If $f$ is a positive constant in $\Omega$, then 
  $H(t)$ is $\beta$\,-concave in $(0,\infty)$ with $\beta=1/(n+2)$. 
\end{itemize}
\end{corollary}
Corollary~\ref{Corollary:1.1} follows from Theorem~\ref{Theorem:1.2} and 
the Borell-Brascamp-Lieb inequality (see Section~5). 
\medskip

We obtain parabolic power concavity properties of the solution $u$ of problem~\eqref{eq:1.1} and \eqref{eq:1.2} 
by developing the method introduced in \cite{cuoghis} and \cite{BS}, 
where quasi-concavity and power concavity properties for elliptic boundary value problems in convex rings and in convex sets were discussed. 
We define $u_{\alpha,p}$ as the \emph{$\alpha$\,-parabolically $p$\,-concave envelope} of the solution~$u$, 
which is the smallest $\alpha$\,-parabolically $p$\,-concave function greater than or equal to $u$, 
and prove that $u_{\alpha,p}=u$ in $\overline{D}$ 
with the aid of the comparison principle for viscosity solutions to \eqref{eq:1.1}. 
This implies that $u$ is $\alpha$\,-parabolically $p$\,-concave in $\overline{D}$. 
Our approach does not require the convexity maximum principle 
and it is completely different from that of \cite{GK}, \cite{IS3}, \cite{IS4}, \cite{Kenni}, and \cite{Korev}.  
\medskip

The rest of this paper is organized as follows. 
In Section~2 we introduce some notation, 
and recall some properties of concave functions and the notion of viscosity solutions. 
In Section~3 we define the $\alpha$\,-parabolically $p$\,-concave envelope 
for nonnegative functions in $D$, then we state and prove the main result of this paper. 
In Section~4 we study the concavity properties of the heat energy $H(t)$ 
with the aid of the Borell-Brascamp-Lieb inequality. 
In Section~5 we apply the results in Sections~3 and 4 to particular parabolic boundary problems, 
and discuss the optimality of our theorems. 
\section{Preliminaries}
In this section we introduce some notation, 
and state some properties of $\alpha$\,-parabolically $p$\,-concave functions. 
Furthermore, we recall the notion of viscosity solutions to \eqref{eq:1.1}. 
\vspace{3pt}

For $x\in{\bf R}^n$ and $r>0$, 
let $B(x,r):=\left\{z\in {\bf R}^n\,:\,|z-x|<r\right\}$. 
For $x=(x_1,\dots,x_n)$, $y=(y_1,\dots,y_n)\in{\bf R}^n$, 
we denote by $x\otimes y$ the $n\times n$ matrix  with entries $(x_iy_j)$, $i,j=1,\dots,n$.
For $m\in\{2,3,\dots\}$, 
we put
$$
\Lambda_m:=\biggr\{\lambda=(\lambda_1,\dots,\lambda_{m})\in(0,1)^m\,:\,
\sum_{i=1}^{m}\lambda_i=1\biggr\}\,.
$$
For $a=(a_1,\dots,a_m)\in(0,\infty)^m$, $\lambda \in\Lambda_m$, and $p\in [-\infty,\,+\infty]$, 
we put 
\begin{equation}
\label{eq:2.1}
 {\bf M}_p(a;\lambda):=
\left\{
\begin{array}{ll}
\left[\lambda_1a_1^p+\lambda_2a_2^p+\cdots+\lambda_ma_m^p\right]^{1/p}
& \mbox{if}\quad p\neq -\infty,\,0,\,+\infty,\vspace{5pt}\\ 
\max\{a_1,\dots,a_m\} & \mbox{if}\quad p=+\infty,\vspace{5pt}\\
a_1^{\lambda_1}\cdots a_m^{\lambda_m} & \mbox{if}\quad p=0,\vspace{5pt}\\
\min\{a_1,a_2,\dots,a_m\} & \mbox{if}\quad p=-\infty,
\end{array}
\right.
\end{equation}
which is the ($\lambda$-weighted) $p$\,-mean of $a$.
For $a=(a_1,\dots,a_m)\in[0,\infty)^m$, we define ${\bf M}_p(a;\lambda)$ as above if $p\geq0$ and 
${\bf M}_p(a;\lambda)=0$ if $p<0$ and $a_1 a_2\cdots a_m=0$. 
Notice that ${\bf M}_p(a;\lambda)$ is a generalization of $M_p(a,b;\lambda)$ defined in Section~1. 
Due to the Jensen inequality, we have 
\begin{equation}
\label{eq:2.2}
{\bf M}_p(a;\lambda)\leq {\bf M}_q(a;\lambda)\quad \mbox{if}\quad -\infty\le p\le q\le\infty,
\end{equation} 
for any $a\in[0,\infty)^m$ and $\lambda \in \Lambda_m$. 
Moreover, it easily follows that
$$
\lim_{p\rightarrow +\infty}
{\bf M}_p(a;\lambda)=\max\{a_1,\dots,a_m\},
\qquad
\lim_{p\rightarrow -\infty}
{\bf M}_p(a;\lambda)=\min\{a_1,\dots,a_m\}.
$$
For further details, see e.g. \cite{HLP}. 
\vspace{5pt}

We state some properties of $\alpha$\,-parabolically $p$\,-concave functions. 
Let $K$ be a convex set in ${\bf R}^n$ , $Q:=K\times(0,\infty)$, $-\infty\le p\le\infty$ and $\alpha\in{\bf R}$. 
The following facts immediately follow from Definition~\ref{Definition:1.2}: 
\begin{itemize}
  \item[(a)] 
  If $v$ is $1$\,-parabolically $p$\,-concave in $K$, 
  then $v$ is $p$\,-concave in $Q$ (in the sense of Definition~\ref{Definition:1.1} with $m=n+1$); 
  \item[(b)] 
  If $w$ is $p$\,-concave in $K$, 
  the function $\tilde{w}(x,t):=w(x)$ is $\alpha$\,-parabolically $p$\,-concave in $Q$ for any $\alpha\in{\bf R}$;
  \item[(c)]   
  If $v$ is $\alpha$\,-parabolically $p$\,-concave in $K$, 
  then $v(\cdot,t)$ is $p$\,-concave in $K$ for any fixed $t>0$.  
\end{itemize}
Furthermore, thanks to \eqref{eq:2.2}, we see that
if $v$ is $\alpha$\,-parabolically $p$\,-concave in $Q$, 
then  
\begin{itemize}
  \item[(d)] $v$ is $\alpha$\,-parabolically $q$\,-concave in $Q$ for any $-\infty\leq q\le p$; 
  \item[(e)] $v$ is $\beta$\,-parabolically $p$\,-concave in $Q$ for any $\beta\ge\alpha$ provided that 
  $u$ is non-decreasing with respect to the time variable~$t$. 
\end{itemize}
In addition, similarly to \cite[Section~2]{Kenni}, 
we see the following: 
\begin{itemize}
  \item[(f)] 
  Let $\{v_j\}$ be nonnegative functions in $Q$ such that for each $j$,  
  $v_j$ is $\alpha_j$\,-parabolically $p_j$\,-concave in $Q$ for some $\alpha_j\in{\bf R}$ and $p_i\in[-\infty,\infty]$. 
  If $v$ is a pointwise limit of a sequence $\{v_j\}$ in $Q$, 
  $\lim_{j\to\infty}\alpha_j=\alpha\in{\bf R}$, and $\lim_{j\to\infty}p_j=p\in[-\infty,\infty]$, 
  then $v$ is $\alpha$\,-parabolically $p$\,-concave in $Q$;
  \item[(g)] Let $\alpha\in{\bf R}$ and $p$, $q\in[0,\infty]$. If $v$ and $w$ are $\alpha$\,-parabolically $p$\,-concave and $q$\,-concave in $Q$, respectively, then $v\cdot w$ is $\alpha$\,-parabolically $r$\,-concave in $Q$, 
  where 
  $$
  \frac{1}{r}=\frac{1}{p}+\frac{1}{q}. 
  $$ 
\end{itemize}

We recall the notion of viscosity 
subsolutions, supersolutions and solutions of \eqref{eq:1.1}. 
Let $\Omega$ be a domain in ${\bf R}^n$, $D=\Omega\times(0,\infty)$, 
and $f=f(x,t,v,\theta)$ a continuous function on $D\times{\bf R}\times {\bf R}^n$. 
An upper semicontinuous function $u$ in $D$ 
is said to be a \emph{viscosity subsolution} of \eqref{eq:1.1} 
if, for any $(x_0,t_0)\in D$, 
the inequality 
$$
\partial_t\phi\le \Delta\phi+f(x,t,\phi,\nabla\phi) 
$$
holds at $(x,t)=(x_0,t_0)$ for all $C^{2,1}(D)$-function $\phi$ satisfying 
$$
\mbox{$\phi(x_0,t_0)=u(x_0,t_0)$ and $\phi(x,t)\ge u(x,t)$ in a neighborhood of $(x_0,t_0)$}.
$$
Analogously, a lower semicontinuous function $u$ in $D$ 
is said to be a \emph{viscosity supersolution} of \eqref{eq:1.1} 
if, for any $(x_0,t_0)\in D$, 
the inequality 
$$
\partial_t\phi\ge \Delta\phi+f(x,t,\phi,\nabla\phi) 
$$
holds at $(x,t)=(x_0,t_0)$ for all $C^{2,1}$ function $\phi$ satisfying 
$$
\mbox{$\phi(x_0,t_0)=u(x_0,t_0)$ and $\phi(x,t)\le u(x,t)$ in a neighborhood of $(x_0,t_0)$.}
$$
A continuous function $u$ in $D$ is said to be a \emph{viscosity solution} of \eqref{eq:1.1} 
if $u$ is a viscosity subsolution and a viscosity supersolution of \eqref{eq:1.1} at the same time.
 
The technique proposed in this paper uses the following (weak) comparison principle for viscosity solutions: 
\begin{equation}
\label{eq:WCP}
\tag{WCP}
\left\{
\begin{array}{l}
\mbox{Let $u\in C(\overline{D})\cap C^{2,1}(D)$ and
$v\in C(\overline{D})$ be a nonnegative classical solution}\vspace{3pt}\\
\mbox{and a nonnegative viscosity subsolution of \eqref{eq:1.1}, respectively, such that}\vspace{3pt}\\
\mbox{$u\geq v$ on $\partial D$. Then $u\geq v $ in $\overline{D}$.}
\end{array}
\right.
\end{equation}
For sufficient conditions for (WCP), see e.g. \cite[Section~8]{UG}. 
\section{Main theorem}
In this section we state and prove the main theorem of this paper, 
which gives a sufficient condition for the solution of \eqref{eq:1.1} and \eqref{eq:1.2} 
to be $\alpha$\,-parabolically $p$\,-concave in $\overline{D}$.  
\vspace{5pt}

We introduce the notion of $\alpha$\,-parabolically $p$\,-concave $\lambda$\,-envelope for nonnegative functions. 
\begin{definition}
\label{Definition:3.1} 
Let $\Omega$ be a bounded convex domain in ${\bf R}^n$, $D:=\Omega\times(0,\infty)$, 
and $u$ a nonnegative function in $\overline{D}$. 
Let $\alpha,p\in[-\infty,+\infty]$ and $\lambda=(\lambda_1,\dots,\lambda_{n+1})\in\Lambda_{n+1}$. 
Then, for any $(x,t)\in\overline{D}$, 
we set 
\begin{eqnarray}
 & & u_{\alpha,p,\lambda}(x,t):=\sup\,
\biggr\{{\bf M}_p(u(y_1,\tau_1),\dots,u(y_{n+1},\tau_{n+1});\lambda)\notag\\
 & & \qquad\qquad
:\, \{(y_i,\tau_i)\}_{i=1}^{n+1}\subset\overline{D},\,\,x=\sum_{i=1}^{n+1}\lambda_i y_i,\,\,
t={\bf M}_\alpha((\tau_1,\dots,\tau_{n+1});\lambda)\biggr\}
\label{eq:3.1}
\end{eqnarray}
 and we call it the $\alpha$\,-parabolically $p$\,-concave $\lambda$\,-envelope of $u$. 
\end{definition}
Taking $y_i=x$ and $\tau_i=t$ for $i=1,\dots,n+1$ in \eqref{eq:3.1}, we easily see that 
\begin{equation}
u_{\alpha,p,\lambda}(x,t)\ge u(x,t),\qquad
(x,t)\in\overline{D}. 
\label{eq:3.2}
\end{equation}
Furthermore, it follows that $u$ is $\alpha$\,-parabolically $p$\,-concave in $\overline{D}$ if and only if 
\begin{equation}
u_{\alpha,p,\lambda}(x,t)\le u(x,t),\qquad
(x,t)\in\overline{D}
\label{eq:3.3}
\end{equation}
(then, equivalently, if and only if $u_{\alpha,p,\lambda}=u$ in $\overline{D}$) for every $\lambda\in\Lambda_{n+1}$. 
We also set
$$
u_{\alpha,p}(x,t)=\sup_{\lambda\in\Lambda_{n+1}}u_{\alpha,p,\lambda}(x,t),
$$
which is the smallest $\alpha$\,-parabolically $p$\,-concave function greater than or equal to $u$ and it is called 
{\em the $\alpha$\,-parabolically $p$\,-concave envelope} of $u$. 
Obviously, $u$ is $\alpha$\,-parabolically $p$\,-concave if and only if
it coincides with $u_{\alpha,p}$.
\vspace{3pt}

We are now ready to state the main theorem of this paper. 
\begin{theorem}
\label{Theorem:3.1}
Let $\Omega$ be a bounded convex domain in ${\bf R}^n$, $D:=\Omega\times(0,\infty)$, and 
$f=f(x,t,v,\theta)$ a nonnegative continuous function in $D\times{\bf R}\times{\bf R}^n$, 
$0<p<1$, and $1/2\le\alpha\le 1$. 
Let $u\in C^{2,1}(D)\cap C(\overline{D})$ satisfy \eqref{eq:1.1} and \eqref{eq:1.2}. 
Assume the following conditions: 
\begin{itemize}
  \item[{\rm (i)}] $\partial_t u\ge 0$ in $D$;
  \item[{\rm (ii)}] For $x\in\overline{\Omega}$ and $y\in\Omega$, 
  \begin{equation}
  \label{eq:3.4}
  \lim_{\rho\to 0+}\frac{\omega_{p,\alpha}(\rho:x,y)}{\rho}=\infty,
  \end{equation}
  where 
  $$
  \omega_{p,\alpha}(\rho:x,y):=u\left(x+\nu\rho,\rho^{1/\alpha}\right)^p,
  \qquad
  \nu:=
  \left\{
  \begin{array}{ll}
  \displaystyle{\frac{y-x}{|y-x|}} & \mbox{if}\quad x\not=y,\vspace{5pt}\\
  0 & \mbox{if}\quad x=y;
  \end{array}
  \right.
  $$
  \item[{\rm (iii)}] The function 
  $$
  g_{\alpha,p,\theta}(x,t,v):=v^{3-1/p}f(x,t^{1/\alpha},v^{1/p},v^{1/p-1}\theta)
  $$
  is concave with respect to $(x,t,v)\in D\times(0,\infty)$ for any fixed $\theta\in{\bf R}^n$. 
\end{itemize}  
Then, for any $\lambda\in\Lambda_{n+1}$, 
the $\alpha$\,-parabolically $p$\,-concave $\lambda$\,-envelope $u_{\alpha,p,\lambda}$ of $u$ is a viscosity subsolution of \eqref{eq:1.1} 
such that $u_{\alpha,p,\lambda}=0$ on $\partial D$. 
In addition, if the comparison principle {\rm (WCP)} holds for equation~\eqref{eq:1.1}, 
then 
$$
u=u_{\alpha,p}\quad\mbox{in}\quad\overline{D}
$$
and $u$ is $\alpha$\,-parabolically $p$\,-concave in $\overline{D}$.
\end{theorem}
We remark that
the function $\omega_{p,\alpha}$ is well-defined for all sufficiently small $\rho>0$ 
since $\Omega$ is convex. 

In order to prove Theorem \ref{Theorem:3.1} we prepare the following two lemmas. 
\begin{lemma}
\label{Lemma:3.1}
Assume the same conditions and notation as in  Theorem~{\rm\ref{Theorem:3.1}}. 
Then
\begin{equation}
\label{eq:3.5}
u(x,t)>0\qquad\mbox{in}\quad D. 
\end{equation}
Furthermore, for any $\lambda\in\Lambda_{n+1}$, 
$u_{\alpha,p,\lambda}$ is continuous on $\overline{D}$ and satisfies  
\begin{equation}
\label{eq:3.6}
u_{\alpha,p,\lambda}>0\qquad\mbox{in}\quad D
\quad\mbox{and}\qquad
u_{\alpha,p,\lambda}=0\quad\mbox{on}\quad \partial D.
\end{equation}
\end{lemma}
{\bf Proof.}
We apply the strong maximum principle 
with the aid of \eqref{eq:3.4} and the nonnegativity of $f$, 
and obtain \eqref{eq:3.5}. 
Moreover, 
by a similar argument to the proof of \cite[Lemma~4.1]{cuoghis} and \cite[Lemma~1]{BS}, 
we see that $u_{\alpha,p,\lambda}\in C(\overline{D})$ and \eqref{eq:3.6}.  
$\Box$
\begin{lemma}
\label{Lemma:3.2}
Assume the same conditions as in Theorem~{\rm\ref{Theorem:3.1}}. 
Then, for any $(x_*,t_*)\in D$ and $\lambda=(\lambda_1,\dots,\lambda_{n+1})\in\Lambda_{n+1}$, 
there exist $\{(x_i,t_i)\}_{i=1}^{n+1}\subset D$ such that 
\begin{eqnarray}
\label{eq:3.7}
 & & x_*=\displaystyle{\sum_{i=1}^{n+1}}\lambda_i x_i,\qquad 
t_*^\alpha=\displaystyle{\sum_{i=1}^{n+1}}\lambda_i t_i^\alpha,\qquad
u_{\alpha,p,\lambda}(x_*,t_*)^p
=\displaystyle{\sum_{i=1}^{n+1}}
\lambda_i u(x_i,t_i)^p,\\
\label{eq:3.8} 
 & & u(x_1,t_1)^{p-1}
\nabla u(x_1,t_1)=\cdots=u(x_{n+1},t_{n+1})^{p-1}\nabla u(x_{n+1},t_{n+1}),\\
\label{eq:3.9}
 & & t_1^{1-\alpha}u(x_1,t_1)^{p-1}\partial_tu(x_1,t_1)=\cdots
=t_{n+1}^{1-\alpha}u(x_{n+1},t_{n+1})^{p-1}\partial_tu(x_{n+1},t_{n+1}).
\end{eqnarray}
\end{lemma}
{\bf Proof.}
Let $(x_*,t_*)\in D$. It follows from \eqref{eq:3.5} that 
\begin{equation}
\label{eq:3.10}
u(x_*,t_*)>0\qquad\mbox{in}\quad D. 
\end{equation}
Since 
$$
\biggr\{(y_1,s_1,y_2,s_2,\dots,y_{n+1},s_{n+1})\in \overline{D}^{n+1}
\,:\,x_*=\sum_{i=1}^{n+1}\lambda_i y_i,\,\,t_*^\alpha=\sum_{i=1}^{n+1}\lambda_i s_i^\alpha\biggr\}
$$
is closed and bounded in $\overline{D}^{n+1}$ and $M_p(u(y_1,s_1),\dots,u(y_{n+1},s_{n+1});\lambda)$ is continuous with respect to $(y_1,s_1y_2,s_2,\dots,y_{n+1},s_{n+1})$, 
we can find $\{(x_i,t_i)\}_{i=1}^{n+1}\subset\overline{D}$ satisfying \eqref{eq:3.7}. 
Since the Lagrange multiplier theorem implies \eqref{eq:3.8}  and \eqref{eq:3.9} 
provided that $\{(x_i,t_i)\}_{i=1}^{n+1}\subset D$, 
it suffices to prove that $\{(x_i,t_i)\}_{i=1}^{n+1}\subset D$. 

The proof is by contradiction. 
Assume that $\{(x_i,t_i)\}_{i=1}^{n+1}\not\subset D$. 
If $\{(x_i,t_i)\}_{i=1}^{n+1}\subset\partial D$, then, by \eqref{eq:1.2}, \eqref{eq:3.2}, and \eqref{eq:3.7} 
we have 
$$
0=u_{\alpha,p,\lambda}(x_*,t_*)\ge u(x_*,t_*),
$$ 
which contradicts \eqref{eq:3.10}. 
This means that $(x_i,t_i)\in\partial D$ and $(x_j,t_j)\in D$ for some $i,j\in\{1,\dots,n+1\}$. 
Here we can assume, without loss of generality, that $i=1$ and $j=2$, that is, 
\begin{equation}
\label{eq:3.11}
(x_1,t_1)\in\partial D,\qquad\qquad (x_2,t_2)\in D. 
\end{equation}
Put 
$$
v(x,\tau):=u(x,\tau^{1/\alpha})^p\quad\mbox{for}\quad (x,\tau)\in\overline{D},
\quad
\tau_*=t_*^\alpha,\quad
\tau_i=t_i^\alpha\quad\mbox{for}\quad i=1,\dots,n+1.
$$  
It follows from \eqref{eq:3.7} that 
\begin{equation}
\label{eq:3.12}
u_{\alpha,p,\lambda}(x_*,t_*)^p=\sum_{i=1}^{n+1}\lambda_i v(x_i,\tau_i),
\qquad
x_*=\sum_{i=1}^{n+1}\lambda_i x_i,
\qquad
\tau_*=\sum_{i=1}^{n+1}\lambda_i\tau_i. 
\end{equation}
These together with the definition of $u_{\alpha,p,\lambda}$ imply that 
\begin{equation}
\label{eq:3.13}
u_{\alpha,p,\lambda}(x_*,t_*)^p\ge \sum_{i=1}^{n+1}\lambda_i v(y_i,\eta_i)
\end{equation}
for all $\{(y_i,\eta_i)\}_{i=1}^{n+1}\subset \overline{D}$ satisfying 
$x_*=\displaystyle{\sum_{i=1}^{n+1}}\lambda_i y_i$ and $\tau_*=\displaystyle{\sum_{i=1}^{n+1}}\lambda_i\eta_i$. 
\vspace{3pt}

Let $\nu:=(x_2-x_1)/|x_2-x_1|$ if $x_1\not=x_2$ and $\nu=0$ if $x_1=x_2$. 
Let $\rho\in(0,1)$, and put 
\begin{equation}
\label{eq:3.14}
\begin{array}{lll}
\tilde{x}_1:=x_1+\nu\displaystyle{\frac{\rho}{\lambda_1}}, 
 & 
\tilde{x}_2:=x_2-\nu\displaystyle{\frac{\rho}{\lambda_2}},
 & 
\tilde{x}_i=x_i\,\,\,(i=3,\dots,n+1),\vspace{5pt}\\ 
\tilde{\tau}_1:=\tau_1+\displaystyle{\frac{\rho}{\lambda_1}}, 
 & 
\tilde{\tau}_2:=\tau_2-\displaystyle{\frac{\rho}{\lambda_2}},
 & 
\tilde{\tau}_i:=\tau_i\quad(i=3,\dots,n+1). 
\end{array}
\end{equation}
It follows from \eqref{eq:3.12} that 
\begin{equation}
\label{eq:3.15}
\sum_{i=1}^{n+1}\lambda_i\tilde{x}_i=\sum_{i=1}^{n+1}\lambda_ix_i=x_*,
\qquad
\sum_{i=1}^{n+1}\lambda_i\tilde{\tau}_i=\sum_{i=1}^{n+1}\lambda_i\tau_i=\tau_*. 
\end{equation}
Furthermore, due to the convexity of $\Omega$, 
we can take a sufficiently small $\rho>0$ so that 
\begin{equation}
\label{eq:3.16}
\{(\tilde{x}_i,\tilde{\tau}_i)\}_{i=1}^{n+1}\subset\overline{D}.  
\end{equation}
Since $u(x_2,t_2)>0$ by \eqref{eq:3.5} and \eqref{eq:3.11}, 
we can find positive constants $M$ and $R_1$ such that 
$$
|(\nabla v)(x,\tau)|+|(\partial_t v)(x,\tau)|\le M\quad\mbox{in}\quad
B(x_2,R_1)\times(\tau_2-R_1,\tau_2+R_1)\subset D. 
$$
Then, taking a sufficiently small $\rho$ if necessary and applying the mean value theorem, 
we obtain 
\begin{equation}
\label{eq:3.17}
\lambda_2[v(\tilde{x}_2,\tilde{\tau}_2)-v(x_2,\tau_2)]
\ge -\lambda_2M|(\tilde{x}_2,\tilde{\tau}_2)-(x_2,\tau_2)|
\ge -2M\rho. 
\end{equation}

On the other hand, 
since $(x_1,t_1)\in\partial D$ by \eqref{eq:3.11}, 
we see that either
$$
{\rm (i)}\quad (x_1,t_1)\in\partial\Omega\times(0,\infty)
\qquad\mbox{or}
\qquad
{\rm (ii)}\quad (x_1,t_1)\in\overline{\Omega}\times\{0\}. 
$$
Consider the case (i). 
Then $x_1\not=x_2$ and $\nu\not=0$. 
The Hopf lemma implies
$$
\liminf_{\rho\to0^+}\frac{u(x_1+\rho\nu,t)}{\rho}>0,\qquad t>0. 
$$
This together with \eqref{eq:1.2} and $p<1$ yields 
$$
(\partial_\nu v)(x,\tilde{\tau}_1)=pu(x,\tilde{\tau}_1^{1/\alpha})^{p-1}(\partial_\nu u)(x,\tilde{\tau}_1^{1/\alpha})
\ge 2(M+1),\quad x\in B(x_1,R_2)\cap D,
$$
for some $R_2>0$. 
Taking a sufficiently small $\rho>0$ if necessary 
and applying the mean value theorem, 
we deduce from \eqref{eq:3.14} that 
\begin{equation}
\label{eq:3.18}
\lambda_1[v(\tilde{x}_1,\tilde{\tau}_1)-v(x_1,\tau_1)]
=\lambda_1 v(\tilde{x}_1,\tilde{\tau}_1)
=\lambda_1[v(\tilde{x}_1,\tilde{\tau}_1)-v(x_1,\tilde{\tau}_1)]
\ge 2(M+1)\rho. 
\end{equation}
Thus, 
by \eqref{eq:3.12}, \eqref{eq:3.17}, and \eqref{eq:3.18} 
we see that 
\begin{eqnarray*}
 & & \sum_{i=1}^{n+1}\lambda_iv(\tilde{x}_i,\tilde{\tau}_i)\\
 & & >[\lambda_1 v(x_1,\tau_1)+2(M+1)\rho]
+[\lambda_2 v(x_2,\tau_2)-2M\rho]
+\sum_{i=3}^{n+1}\lambda_iv(x_i,\tau_i)\\
& & >\sum_{i=1}^{n+1}\lambda_iv(x_i,\tau_i)=u_{\alpha,p,\lambda}(x_*,t_*)^p,
\end{eqnarray*}
which together with \eqref{eq:3.15} contradicts \eqref{eq:3.13}. 

Consider the case~(ii). 
By \eqref{eq:3.4} and \eqref{eq:3.14} 
we see that 
\begin{equation}
\lambda_1[v(\tilde{x}_1,\tilde{\tau}_1)-v(x_1,t_1)]
=\lambda_1v(\tilde{x}_1,\tilde{\tau}_1)
=\lambda_1\omega_{p,\alpha}(\lambda_1^{-1}\rho:x_1,x_2)
>2M\rho
\label{eq:3.19}
\end{equation}
for all sufficiently small $\rho$. 
Therefore, 
taking a sufficiently small $\rho>0$ if necessary and 
combining \eqref{eq:3.19} with \eqref{eq:3.17}, 
we have 
\begin{eqnarray}
 & & \sum_{i=1}^{n+1}\lambda_iv(\tilde{x}_i,\tilde{\tau}_i)\notag\\
 & & >\lambda_1v(x_1,\tau_1)+2M\rho
+\lambda_2v(x_2,\tau_2)-2M\rho
+\sum_{i=3}^{n+1}\lambda_iv(x_i,\tau_i)\notag\\
 & & =\sum_{i=1}^{n+1}\lambda_iv(x_i,\tau_i)=u_{\alpha,p,\lambda}(x_*,t_*)^p.
\end{eqnarray}
This together with \eqref{eq:3.15} contradicts \eqref{eq:3.13}. 
Therefore, in the both cases (i) and (ii), 
we have a contradiction. Thus we see that $\{(x_i,t_i)\}_{i=1}^{n+1}\subset D$, 
and Lemma~\ref{Lemma:3.2} follows. 
$\Box$
\vspace{5pt}

Now we are ready to prove Theorem~\ref{Theorem:3.1}. 
\vspace{5pt}
\newline
{\bf Proof of Theorem~\ref{Theorem:3.1}.} 
Let $(x_*,t_*) \in D$ and $\lambda=(\lambda_1,\dots,\lambda_{n+1})\in\Lambda_{n+1}$. 
By Lemma~\ref{Lemma:3.2} 
we can find $\{(x_i,t_i)\}_{i=1}^{n+1}\subset D$ satisfying \eqref{eq:3.7}.  
Put 
\begin{equation}
\label{eq:3.21}
\begin{array}{l}
v_*:=u_{\alpha,p,\lambda}(x_*,t_*)^p,
\qquad
v_i:=u(x_i,t_i)^p,
\qquad
\displaystyle{a_i:=\frac{v_i}{v_*}},\vspace{5pt}\\
y_i(x)=x_i+a_i(x-x_*),
\qquad
\tau_i(t)=\left[t_i^\alpha+a_i(t^\alpha-t_*^\alpha)\right]^{1/\alpha},
\vspace{5pt}
\end{array}
\end{equation}
for $x\in{\bf R}^n$, $t\ge 0$, and $i=1,\dots,n+1$. 
These imply that 
\begin{equation}
\label{eq:3.22}
v_*=\sum_{i=1}^{n+1}\lambda_i v_i,\quad
\sum_{i=1}^{n+1}\lambda_i a_i=1,\quad
x=\sum_{i=1}^{n+1}\lambda_iy_i(x),\quad
t^\alpha=\sum_{i=1}^{n+1}\lambda_i\tau_i(t)^\alpha.
\end{equation}
Furthermore, we see that the function 
\begin{equation}
\label{eq:3.23}
\varphi(x,t):=
\biggr[\,\sum_{i=1}^{n+1}\lambda_iu(y_i(x),\tau_i(t))^p\,\biggr]^{1/p}
\end{equation}
is a $C^{2,1}$-function in a neighborhood of $(x_*,t_*)\in D$ and satisfies 
\begin{equation}
\varphi(x_*,t_*)^p=\sum_{i=1}^{n+1}\lambda_i u(x_i,t_i)^p
=u_{\alpha,p,\lambda}(x_*,t_*)^p =v_*.\vspace{5pt}
 \label{eq:3.24}
\end{equation}
Moreover, it follows from the definition of $u_{\alpha,p,\lambda}$ and \eqref{eq:3.22} that  
\begin{equation}
u_{\alpha,p,\lambda}(x,t)\ge \varphi(x,t)
\label{eq:3.25}
\end{equation}
in a neighborhood of $(x_*,t_*)$. 

We prove 
\begin{equation}
\partial_t\varphi(x_*,t_*)
\leq\Delta\varphi(x_*,t_*)+
f(x_*,t_*,\varphi(x_*,t_*),\nabla\varphi(x_*,t_*)).
\label{eq:3.26}
\end{equation}
By \eqref{eq:3.23} we have 
\begin{eqnarray}
 & & \nabla\varphi(x,t)
=\varphi(x,t)^{1-p}
\sum_{i=1}^{n+1}\lambda_ia_i
u(y_i(x),\tau_i(t))^{p-1}\nabla u(y_i(x),\tau_i(t)),\label{eq:3.27}\\
 & & \nabla^2\varphi(x,t)
=(1-p)\varphi(x,t)^{-1}\nabla\varphi(x,t)\otimes\nabla\varphi(x,t)\notag\\
 & & \qquad
 -(1-p)\varphi(x,t)^{1-p}
 \sum_{i=1}^{n+1}\lambda_ia_i^2
u(y_i(x),\tau_i(t))^{p-2}\nabla u(y_i(x),\tau_i(t))\otimes\nabla u(y_i(x),\tau_i(t))\notag\\
 & & \qquad
+\varphi(x,t)^{1-p}\sum_{i=1}^{n+1}\lambda_ia_i^2
u(y_i(x),\tau_i(t))^{p-1}\nabla^2 u(y_i(x),\tau_i(t)),\label{eq:3.28}
\end{eqnarray}
in a neighborhood of $(x_*,t_*)$. 
Since $y_i(x_*)=x_i$ and $\tau_i(t_*)=t_i$, 
by \eqref{eq:3.8}, \eqref{eq:3.22}, \eqref{eq:3.24}, and \eqref{eq:3.27} 
we have 
\begin{eqnarray}
\nabla\varphi(x_*,t_*)
\!\!\! & = &\!\!\! 
\varphi(x_*,t_*)^{1-p}\sum_{i=1}^{n+1}\lambda_i a_iu(x_i,t_i)^{p-1}\nabla u(x_i,t_i)\notag\\
\!\!\! & = &\!\!\!
\varphi(x_*,t_*)^{1-p}u(x_i,t_i)^{p-1}\nabla u(x_i,t_i)\notag\\
\!\!\! & = &\!\!\!
v_*^{1/p-1}v_i^{-1/p+1}\nabla u(x_i,t_i)
\label{eq:3.29}
\end{eqnarray}
for $i=1,\dots,n+1$. 
This together with \eqref{eq:3.21} implies that 
\begin{eqnarray}
f(x_i,t_i,u(x_i,t_i),\nabla u(x_i,t_i))
 \!\!\! & = & \!\!\! 
 f(x_i,t_i,v_i^{1/p},v_*^{1-1/p}\nabla\varphi(x_*,t_*)v_i^{1/p-1})\notag\\
 \!\!\! & = & \!\!\! 
 v_i^{-3+1/p}g_{\alpha,p,\theta}(x_i,t_i^\alpha,v_i)
 \label{eq:3.30}
\end{eqnarray}
with $\theta:=v_*^{1-1/p}\nabla\varphi(x_*,t_*)$, 
where $i=1,\dots,n+1$. 
Furthermore, by \eqref{eq:3.21} and \eqref{eq:3.28} we obtain  
\begin{eqnarray*}
& & \nabla^2\varphi(x_*,t_*)
=(1-p)\varphi(x_*,t_*)^{-1}\nabla\varphi(x_*,t_*)\otimes\nabla\varphi(x_*,t_*)\notag\\
 & & \qquad
 -(1-p)\varphi(x_*,t_*)^{-1+p}
 \sum_{i=1}^{n+1}\lambda_i\biggr(\frac{v_i}{v_*}\biggr)^2
v_i^{1-2/p}\nabla u(x_i,t_i)\otimes\nabla u(x_i,t_i)\notag\\
 & & \qquad\qquad
+\varphi(x_*,t_*)^{1-p}\sum_{i=1}^{n+1}\lambda_i\biggr(\frac{v_i}{v_*}\biggr)^2
v_i^{1-1/p}\nabla^2 u(x_i,t_i).
\end{eqnarray*}
This together with \eqref{eq:3.22}, \eqref{eq:3.24}, and \eqref{eq:3.29} implies that 
\begin{eqnarray*}
\nabla^2\varphi(x_*,t_*)
\!\!\! & = &\!\!\! 
(1-p)v_*^{-1/p-1}
\biggr[v_*-\sum_{i=1}^{n+1}\lambda_iv_i\biggr]
\nabla\varphi(x_*,t_*)\otimes\nabla\varphi(x_*,t_*)\\
 & &\qquad\qquad\qquad\qquad\quad
 +\sum_{i=1}^{n+1}
\lambda_i\frac{v_i^{3-1/p}}{v_*^{3-1/p}}\nabla^2u(x_i,t_i)\\
\!\!\! & = &\!\!\!
\sum_{i=1}^{n+1}
\lambda_i\frac{v_i^{3-1/p}}{v_*^{3-1/p}}\nabla^2u(x_i,t_i).
\end{eqnarray*}
Then it follows from \eqref{eq:1.1} that 
\begin{eqnarray}
\Delta\varphi(x_*,t_*)
\!\!\! & = &\!\!\!\sum_{i=1}^{n+1}
\lambda_i\frac{v_i^{3-1/p}}{v_*^{3-1/p}}\Delta u(x_i,t_i)\notag\\
\!\!\! & = &\!\!\!
\sum_{i=1}^{n+1}
\lambda_i\frac{v_i^{3-1/p}}{v_*^{3-1/p}}\biggr[\partial_tu(x_i,t_i)-f(x_i,t_i,u(x_i,t_i),\nabla u(x_i,t_i))\biggr]. 
\label{eq:3.31}
\end{eqnarray}
On the other hand, since
\begin{equation}
\label{eq:3.32}
\partial_t\varphi(x_*,t_*)
=\varphi(x_*,t_*)^{1-p}
\sum_{i=1}^{n+1}\lambda_i a_i\biggr(\frac{t_i}{t_*}\biggr)^{1-\alpha}u(x_i,t_i)^{p-1}\partial_tu(x_i,t_i),
\end{equation}
similarly to \eqref{eq:3.29}, 
by \eqref{eq:3.9} and \eqref{eq:3.22} we have 
$$
\partial_t\varphi(x_*,t_*)
=\biggr(\frac{t_i}{t_*}\biggr)^{1-\alpha}\varphi(x_*,t_*)^{1-p}
u(x_i,t_i)^{p-1}\partial_tu(x_i,t_i). 
$$
This together with \eqref{eq:3.21} and \eqref{eq:3.24} yields
\begin{equation}
\label{eq:3.33}
\partial_tu(x_i,t_i)=
\biggr(\frac{t_i^\alpha}{t_*^\alpha}\biggr)^{1-1/\alpha}
v_*^{1-1/p}v_i^{1/p-1}\partial_t\varphi(x_*,t_*),
\qquad i=1,\dots,n+1. 
\end{equation}
Therefore we deduce from \eqref{eq:3.30}, \eqref{eq:3.31}, and \eqref{eq:3.33} that 
\begin{eqnarray}
 & & \partial_t\varphi(x_*,t_*)-\Delta\varphi(x_*,t_*)\notag\\
 & & =\partial_t\varphi(x_*,t_*)
 \left[1-\sum_{i=1}^{n+1}
\lambda_i \frac{v_i^2}{v_*^2}\biggr(\frac{t_i^\alpha}{t_*^\alpha}\biggr)^{1-1/\alpha}\right]
+\sum_{i=1}^{n+1}\lambda_i \frac{g_{\alpha,p,\theta}(x_i,t_i^\alpha,v_i)}{v_*^{3-1/p}}.
\label{eq:3.34}
\end{eqnarray}
On the other hand, since $1/2\le\alpha\le 1$, 
$h(\eta,\tau):=\eta^2 \tau^{1-1/\alpha}$ is a convex function in ${\bf R}\times(0,\infty)$ 
(see e.g. \cite[Lemma~A.1~(i)]{IS3}). 
Then, due to \eqref{eq:3.22}, 
we have 
\begin{eqnarray}
\sum_{i=1}^{n+1}
\lambda_i \frac{v_i^2}{v_*^2}\biggr(\frac{t_i^\alpha}{t_*^\alpha}\biggr)^{1-1/\alpha}
\!\!\! & = &\!\!\!
\sum_{i=1}^{n+1}
\lambda_i h\biggr(\frac{v_i}{v_*},\frac{t_i^\alpha}{t_*^\alpha}\biggr)\notag\\
\!\!\! & \ge &\!\!\!
h\biggr(\sum_{i=1}^{n+1}
\lambda_i\frac{v}{v_*}_i,\sum_{i=1}^{n+1}
\lambda_i\frac{t_i^\alpha}{t_*^\alpha}\biggr)
=h(1,1)=1.
\label{eq:3.35}
\end{eqnarray}
Furthermore, 
since $g_{\alpha,p,\theta}=g_{\alpha,p,\theta}(x,t,v)$ is concave in $D\times(0,\infty)$, 
by \eqref{eq:3.22} we have 
\begin{eqnarray}
\sum_{i=1}^{n+1}\lambda_i \frac{g_{\alpha,p,\theta}(x_i,t_i^\alpha,v_i)}{v_*^{3-1/p}}
\!\!\! & \le &\!\!\! 
\frac{1}{v_*^{3-1/p}}g_{\alpha,p,\theta}
\biggr(\,\sum_{i=1}^{n+1}\lambda_ix_i,\sum_{i=1}^{n+1}\lambda_it_i^\alpha,\sum_{i=1}^{n+1}\lambda_iv_i,\biggr)
\notag\\
\!\!\! & = &\!\!\!\frac{g_{\alpha,p,\theta}(x_*,t_*^\alpha,v_*)}{v_*^{3-1/p}}
=f(x_*,t_*,v_*,v_*^{1/p-1}\theta)\notag\\
\!\!\! & = &\!\!\!
f(x_*,t_*,\varphi(x_*,t_*),\nabla\varphi(x_*,t_*)). 
\label{eq:3.36}
\end{eqnarray}
Therefore, applying \eqref{eq:3.32}, \eqref{eq:3.35} and \eqref{eq:3.36} 
to \eqref{eq:3.34} (and taking in account assumption~(i) on the sign of $\partial_t u$), we obtain \eqref{eq:3.26}. 

Since $(x_*,t_*)$ is arbitrary, 
by \eqref{eq:3.24}--\eqref{eq:3.26} and Lemma~\ref{Lemma:3.1} 
we see that 
$u_{\alpha,p,\lambda}$ is a viscosity subsolution of \eqref{eq:1.1} such that $u_{\alpha,p,\lambda}=0$ on $\partial D$. 
Furthermore, if the comparison principle (WCP) holds for \eqref{eq:1.1}, 
then we obtain \eqref{eq:3.3} for all $(x,t)\in\overline{D}$ and $\lambda\in\Lambda_{n+1}$. 
This together with \eqref{eq:3.2} implies that 
$$
u=u_{\alpha,p,\lambda}=u_{\alpha,p}\qquad\mbox{in}\qquad\overline{D}.
$$
Hence, $u$ is $\alpha$\,-parabolically $p$\,-concave in $\overline{D}$, 
and Theorem~\ref{Theorem:3.1} follows.
$\Box$
\begin{remark}
\label{Remark:3.1} 
If $f$ is independent of the time variable $t$, 
then condition~{\rm (iii)} in Theorem~{\rm\ref{Theorem:3.1}} coincides with the following:   
\begin{itemize}
  \item[{\rm (iii')}] The function 
  $$
  g_{p,\theta}(x,v):=v^{3-1/p}f(x,v^{1/p},v^{1/p-1}\theta)
  $$
  is concave with respect to $(x,v)\in \Omega\times(0,\infty)$ for any fixed $\theta\in{\bf R}^n$. 
\end{itemize}
This condition has already been used in {\rm\cite[Theorem~3.3]{Kenni}} 
for the study of power concavity properties of the solutions of 
$$
\Delta v+f(x,v,\nabla v)=0\quad\mbox{in}\quad\Omega,
\qquad
v=0\quad\mbox{in}\quad\Omega. 
$$
\end{remark}
\section{Concavity of heat energy}
One can consider problem \eqref{eq:1.1} as a mathematical model describing the following situation: 
a cold convex body $\Omega$ with homogeneous density is immersed in liquid kept at constant zero temperature
and is heated by the source term $f$. 
Then $u(x,t)$ describes the temperature at the point $x\in\Omega$ at time $t$, 
and the quantity
$$
H(t)=\int_\Omega u(x,t)\,dx
$$
represents the heat energy of $\Omega$, up to multiplication by a constant. 
In this section we prove the following theorem on power concavity properties of the heat energy 
for parabolically power concave functions.
\begin{theorem}
\label{Theorem:4.1}
Let $\Omega$ be a bounded convex domain in ${\bf R}^n$, $D:=\Omega\times(0,\infty)$, 
$0<\alpha\leq 1$, and $p\ge -1/n$. 
If $u$ is $\alpha$\,-parabolically $p$\,-concave and 
non-decreasing with respect to the time variable~$t$ in $D$, 
then $H(t)$ is $q$\,-concave with 
\begin{equation}
\label{eq:4.1}
q=\left\{
\begin{array}{ll}
1/n &\text{if}\quad p=+\infty,\vspace{5pt}\\
p/(np+1) &\text{if}\quad p\in(-1/n,+\infty),\vspace{5pt}\\
-\infty&\text{if}\quad p=-1/n.
\end{array}
\right.
\end{equation}
\end{theorem}
In order to prove Theorem~\ref{Theorem:4.1}, 
we recall the Borell-Brascamp-Lieb inequality, which is a generalization of the Pr\'ekopa-Leindler inequality. 
See \cite[Theorem 10.1]{Gardner}. 
\begin{proposition}
\label{Proposition:4.1}
Let $\lambda\in(0,1)$, $f,g,h$ nonnegative functions in $L^1({\bf R}^n)$, 
and $-1/n\leq p\leq \infty$. 
Assume that
\begin{equation}
\label{eq:4.2}
h\big((1-\lambda)x+\lambda y\big)\geq M_p(f(x),g(y);\lambda)
\end{equation}
for all $x\in\text{\em sprt}(f),\,y\in\text{\em sprt}(g)$. Then
$$
\int_{{\bf R}^n} h\,dx\ge M_q\left(\int_{{\bf R}^n} f\,dx,\int_{{\bf R}^n} g\,dx\,;\lambda\right)\,,
$$
where $q$ is as in \eqref{eq:4.1}.
\end{proposition}
The Pr\'ekopa-Leindler inequality corresponds to the case $p=0$. 

Notice that usually, see \cite{BL, Gardner}, 
assumption~\eqref{eq:4.2} is required to hold for every $x,y\in{\bf R}^n$, 
but there the definition of $M_p$ is different in that $M_p(a,b)=0$ 
as soon as $ab=0$ even for $p>0$. 
On the other hand, 
this makes $M_p$ not continuous for $p>0$ 
and here we prefer to work with continuous $p$-means for several reasons.
\vspace{5pt}
\newline
{\bf Proof Theorem~\ref{Theorem:4.1}.}
Since $u$ is $\alpha$\,-parabolically $p$\,-concave and non-deceasing with respect to the time variable, 
it follows from property~(e) in Section~2 that 
$u$ is $1$\,-parabolically $p$\,-concave in $D$. 
This implies that 
for any $\lambda\in(0,1)$ and $r$, $s>0$, 
the inequality 
\begin{equation}
\label{eq:4.3}
u\big((1-\lambda)x+\lambda y,(1-\lambda)r+\lambda s\big)\ge M_p\big(u(x,r),u(y,s);\lambda\big)
\end{equation}
holds for all $x$, $y\in\Omega$. 
Then, by Proposition~\ref{Proposition:4.1} we have 
$$
H((1-\lambda)r+\lambda s)\ge M_q\left(H(r), H(s);\lambda\right),
$$
where $q$ is as in \eqref{eq:4.1}. 
Therefore we see that $H(t)$ is $q$\,-concave in $(0,\infty)$, 
and Theorem~\ref{Theorem:4.1} follows.
$\Box$\vspace{5pt}

\begin{remark}
\label{Remark:4.1} 
Assume the same conditions as in Theorem~{\rm\ref{Theorem:4.1}}. 
Similarly to the proof of Theorem~{\rm\ref{Theorem:4.1}}, 
we have the following: 
\vspace{3pt}
\newline
{\rm (i)} 
$H(t^{1/\alpha})$ is $q$\,-concave without the assumption on the monotonicity of $u$ with respect to the time variable;
\vspace{3pt}
\newline
{\rm (ii)} 
Let $m>0$ and $p\ge -m/n$. It follows from \eqref{eq:4.3} that
$$
u\big((1-\lambda)x+\lambda y,(1-\lambda)r+\lambda s\big)^m\ge M_{p/m}\big(u(x,r)^m,u(y,s)^m;\lambda\big)
$$
for all $x$, $y\in\Omega$, and we see that  
$$
H_m(t)=\left(\int_\Omega u(x,t)^mdx\right)^{1/m},\qquad m>0. 
$$
is $q$\,-concave with 
$$
q=\left\{
\begin{array}{ll}
1/n &\text{if}\quad p=+\infty,\vspace{5pt}\\
p/(np+m) &\text{if}\quad p\in(-m/n,+\infty),\vspace{5pt}\\
-\infty&\text{if}\quad p=-m/n.
\end{array}
\right.
$$
\end{remark}
\section{Applications}
In this section we apply Theorems~\ref{Theorem:3.1} and \ref{Theorem:4.1}
to particular parabolic boundary value problems and discuss the sharpness of our results. 
We first deal with the case where $f(x,t,v,\theta)$ is independent of $v$ and $\theta$, 
and prove the following theorem. 
Theorem~\ref{Theorem:1.2} and Corollary~\ref{Corollary:1.1} 
easily follow from Theorem~\ref{Theorem:5.1} with $\gamma=0$ and property~(c) in Section~2. 
\begin{theorem}
\label{Theorem:5.1}
Let $\Omega$ be a bounded convex domain in ${\bf R}^n$, $D:=\Omega\times(0,\infty)$, 
and $f$ a nonnegative function in $\Omega$. 
Let $u\in C^{2,1}(D)\cap C(\overline{D})$ satisfy 
\begin{equation}
\label{eq:5.1}
\partial_tu=\Delta u +t^\gamma f(x)\quad\mbox{in}\quad D
\qquad\mbox{and}\qquad
u=0\quad\mbox{on}\quad\partial D,
\end{equation}
where $0\le\gamma\le 1/2$. 
\vspace{3pt}
\newline
{\rm (i)} If $f$ is $q$\,-concave in $\Omega$ for some $q\ge 1$, 
then $u$ is parabolically $p$\,-concave in $\overline{D}$ with 
$$
p=\frac{q}{1+2q+2\gamma q}. 
$$
Furthermore, the heat energy $H(t)$ is $r$\,-concave in $(0,\infty)$ with 
$$
r=\frac{q}{(n+2+\gamma)q+1};
$$
{\rm (ii)} If $f$ is a positive constant in $\Omega$, 
then $u$ is parabolically $p$\,-concave in $\overline{D}$ with 
$$
p=\frac{1}{2(1+\gamma)}
$$ 
and $H(t)$ is $r$\,-concave in $(0,\infty)$ with $r=1/(n+2+\gamma)$. 
\end{theorem}
In order to prove Theorem~\ref{Theorem:5.1}, 
we prepare the following lemma.
\begin{lemma}
\label{Lemma:5.1} 
Let $\Omega$ be a bounded convex smooth domain in ${\bf R}^n$, $D:=\Omega\times(0,\infty)$, 
$d\ge 0$, and $\gamma\in[0,1)$. 
Let $U\in C^2(D)\cap C(\overline{D})$ satisfy 
\begin{equation}
\label{eq:5.2}
\partial_t U=\Delta U+t^\gamma\mbox{{\rm dist}}(x,\partial\Omega)^d\quad\mbox{in}\quad D,
\qquad
U=0\quad\mbox{on}\quad\partial D.
\end{equation}
Let $x_*\in\overline{\Omega}$ and $y_*\in\Omega$, and put 
$\nu:=(y_*-x_*)/|y_*-x_*|$ if $y_*\neq x_*$ and $\nu=0$ if $y_*=x_*$. 
Then there exists a constant $C$ such that 
\begin{equation}
\label{eq:5.3}
U(x_*+\nu\rho,\rho^2)\ge C\rho^{2\gamma+d+2}
\end{equation}
for all sufficiently small $\rho>0$.
\end{lemma}
{\bf Proof.}
By the use of the Dirichlet heat kernel $G=G(x,y,t)$ on $\Omega$, 
the function $U$ is represented by 
\begin{equation}
\label{eq:5.4}
U(x,t)=\int_0^t\int_\Omega G(x,y,t-s)s^\gamma\mbox{{\rm dist}}(y,\partial\Omega)^d\,dyds.
\end{equation}
Let $x_*\in\overline{\Omega}$ and $y_*\in\Omega$, with $y_*\neq x_*$. 
Due to the convexity of the domain $\Omega$, 
we can find an open convex cone $K$ in ${\bf R}^n$ with the vertex at the origin 
such that 
\begin{eqnarray*}
 & & \nu=(y_*-x_*)/|y_*-x_*|\in K\cap\,\partial B(0,1),\vspace{3pt}\\
 & & x_*+(K\cap B(0,R))\subset\overline{\Omega},\vspace{3pt}\\
 & & \mbox{{\rm dist}}\,(x,\partial\Omega)\ge C_1|x-x_*|
\quad\mbox{if}\quad x-x_*\in K\cap B(0,R),
\end{eqnarray*}
for some positive constants $R$ and $C_1$. 
These together with \eqref{eq:5.4} imply that 
\begin{eqnarray*}
U(x,t) \!\!\! & \ge &\!\!\!
 C_2t^\gamma\int_{t/4}^{t/2}\int_{x_*+K\cap B(0,R)}
G(x,y,t-s)|y-x_*|^d\,dyds\\
\!\!\! & = &\!\!\!
C_2t^\gamma\int_{t/4}^{t/2}\int_{K\cap B(0,R)}
G(x,x_*+y,t-s)|y|^d\,dyds
\end{eqnarray*}
for some constant $C_2>0$. 
Then we have 
\begin{eqnarray}
 & & U(x_*+\nu\rho,\rho^2)\notag\\
 & & 
 \ge C_2\rho^{2\gamma}
 \int_{\rho^2/4}^{\rho^2/2}\int_{K\cap B(0,R)}
G(x_*+\rho\nu,x_*+y,\rho^2-s)|y|^d\,dyds\notag\\
 & & 
 \ge C_2\rho^{2\gamma+d+2}
 \int_{1/4}^{1/2}\int_{K\cap B(0,2)}
 G_\rho(\nu,z,1-\eta)|z|^d\,dzd\eta
 \label{eq:5.5}
\end{eqnarray}
for all sufficiently small $\rho>0$, where 
$$
G_\rho(x,y,t):=\rho^n G(x_*+\rho x,x_*+\rho y,\rho^2 t). 
$$

Consider the case $x_*\in\Omega$. 
Let $\Gamma$ be the Gauss kernel, that is, 
\begin{equation}
\label{eq:5.6}
\Gamma(x,y,t):=(4\pi t)^{-\frac{n}{2}}\exp\left(-\frac{|x-y|^2}{4t}\right). 
\end{equation}
Since $G_\rho(x,y,t)$ is  the Dirichlet heat kernel in $\Omega_\rho:=\rho^{-1}(\Omega-x_*)$, 
the function 
$$
w(x,t):=\Gamma(x,y,t)-G_\rho(x,y,t). 
$$
satisfies
$$ 
\left\{
\begin{array}{ll}
\partial_t w=\Delta w & \mbox{in}\quad\Omega_\rho\times(0,\infty),\vspace{3pt}\\
w(x,t)=\Gamma(x,y,t) & \mbox{on}\quad\partial\Omega_\rho\times(0,\infty),\vspace{3pt}\\
w(x,0)=0 & \mbox{in}\quad\Omega_\rho. 
\end{array}
\right.
$$
Since $\Omega_\rho$ tends to ${\bf R}^n$ as $\rho\to 0$, 
it follows from the maximum principle and \eqref{eq:5.6} that 
\begin{equation}
\label{eq:5.7}
\lim_{\rho\to 0}\sup_{x\in\Omega_\rho,y\in E,t\in(0,T)}|\Gamma(x,y,t)-G_\mu(x,y,t)|
\le\lim_{\rho\to 0}\sup_{x\in\partial\Omega_\rho,y\in E,t\in(0,T)}\Gamma(x,y,t)=0
\end{equation}
for any compact set $E$ in ${\bf R}^n$ and $T>0$. 
This means that 
$$
\lim_{\rho\to 0}G_\rho(x,y,1-\eta)
=\Gamma(x,y,1-\eta)
$$
uniformly for all $y\in B(0,2)$ and $\eta\in(1/4,1/2)$. 
Therefore, by \eqref{eq:5.5} we can find a positive constant $C_3$ such that 
$$
U(x_*+\nu\rho,\rho^2)\ge C_3\rho^{2\gamma+d+2}
$$
for all sufficiently small $\rho>0$. 
Thus \eqref{eq:5.3} holds in the case $x_*\in\Omega$, $y_*\neq x_*$.  

If $y_*=x_*$ the proof works in the same way (in fact with some simplification).

Next we consider the case $x_*\in\partial\Omega$. 
Due to the regularity of $\Omega$, 
$\Omega_\rho$ tends to an open half space $\Pi$ with $0\in\partial\Pi$ as $\rho\to 0$. 
Let $G_D=G_D(x,y,t)$ be the Dirichlet heat kernel on $\Pi$. 
Then, similarly to the case $x_*\in\Omega$, 
we see that 
$$
\lim_{\rho\to 0}\sup_{x\in\Omega_\rho,y\in E,t\in(0,T)}|\Gamma_D(x,y,t)-G_\mu(x,y,t)|
\le\lim_{\rho\to 0}\sup_{x\in\partial\Omega_\rho,y\in E,t\in(0,T)}\Gamma_D(x,y,t)=0
$$
for any compact set $E$ in $\Pi$ and $T>0$,  
and obtain \eqref{eq:5.3}. 
Thus Lemma~\ref{Lemma:5.1} follows. 
$\Box$\vspace{5pt}
\newline
{\bf Proof of Theorem~\ref{Theorem:5.1}.}
We prove assertion~(i). 
If $f\equiv 0$ in $\Omega$, then $u\equiv 0$ in $D$, and assertion~(i) easily follows.  
So it suffices to consider the case where $f\not\equiv 0$ in $\Omega$. 
Then, by the concavity of $f$ we see that 
\begin{equation}
\label{eq:5.8}
f>0\quad\mbox{in}\quad\Omega.
\end{equation}
On the other hand, we can assume, without loss of generality, 
that $\Omega$ is smooth and $f$ is $q$\,-concave in $\overline{\Omega}$. 
Indeed, there exists a sequence of smooth domains $\{\Omega_k\}_{k=1}^\infty$ such that 
$$
\Omega_1\subset\Omega_2\subset\cdots\subset \Omega_k\subset\cdots
\qquad\mbox{and}\qquad
\bigcup_{k=1}^\infty\Omega_k=\Omega. 
$$
For any $k=1,2,\dots$, since $f$ is bounded and locally Lipschitz in $\overline{\Omega_k}$, 
there exists a classical solution $u_k$ of 
\begin{equation}
\label{eq:5.9}
\partial_tu=\Delta u+t^\gamma f(x)\quad\mbox{in}\quad D_k,
\qquad
u=0\quad\mbox{on}\quad\partial D_k,
\end{equation}
where $D_k:=\Omega_k\times(0,\infty)$, such that 
\begin{equation}
\label{eq:5.10}
u_k\in C(\overline{D_k})\cap C^{2,1}(\overline{D_k}\setminus[\partial\Omega_k\times\{0\}]). 
\end{equation}
By the comparison principle we see that 
$$
0\le u_k(x,t)\le u_{k+1}(x,t)\le u(x,t)\qquad\mbox{in}\quad D_k.
$$
Furthermore, by the regularity theorems for parabolic equations and the uniqueness of the solution of \eqref{eq:1.7} 
we have 
\begin{equation}
\label{eq:5.11}
\lim_{k\to\infty}u_k(x,t)=u(x,t)\qquad\mbox{in}\quad D. 
\end{equation}
This means that 
if $u_k$ is $\alpha$\,-parabolically $p$\,-concave in $\overline{D_k}$ for all $k=1,2,\dots$, 
then $u$ is $\alpha$\,-parabolically $p$\,-concave in $\overline{D}$. 
Therefore it suffices to prove assertion~(i) 
in the case where $\Omega$ is smooth and $f$ is $q$\,-concave in $\overline{\Omega}$. 

Assume then that $\Omega$ is smooth 
and $f$ is $q$\,-concave in $\overline{\Omega}$. 
Similarly to \eqref{eq:5.10}, we see that 
$$
u\in C(\overline{D})\cap C^{2,1}(\overline{D}\setminus[\partial\Omega\times\{0\}]). 
$$
Furthermore, since $\hat{u}:=\partial_t u$ satisfies 
$$
\left\{
\begin{array}{ll}
\partial_t\hat{u}\ge\Delta \hat{u} & \mbox{in}\quad D,\vspace{3pt}\\
\hat{u}=0 & \mbox{on}\quad\partial\Omega\times(0,\infty),\vspace{3pt}\\
\hat{u}(x,0)\ge 0 & \mbox{in}\quad\Omega,
\end{array}
\right.
$$
with a strict inequality in the first equation if $\gamma >0$ or in the last equation if $\gamma=0$,
the comparison principle implies that 
\begin{equation}
\label{eq:5.12}
\partial_t u=\hat{u}>0\qquad\mbox{in}\quad D.
\end{equation}

On the other hand, 
since $f$ is $q$\,-concave in $\overline{\Omega}$, 
by \eqref{eq:5.8} we can find a positive constant $C_1$ such that 
$$
f(x)\ge C_1\mbox{{\rm dist}}(x,\partial\Omega)^{1/q},\qquad x\in\Omega. 
$$
Then it follows from the comparison principle that 
\begin{equation}
\label{eq:5.13}
u(x,t)\ge C_2 U(x,t)\quad\mbox{in}\quad D
\end{equation}
for some positive constant $C_2$, where $U$ is a solution of \eqref{eq:5.2} with $d=1/q$. 
Let 
\begin{equation}
\label{eq:5.14}
0<p<\frac{q}{1+2q+2\gamma q}
\qquad\mbox{and}\qquad
0\le\gamma<\frac{1}{2}.
\end{equation}
By Lemma~\ref{Lemma:5.1} and \eqref{eq:5.13}, 
for any $x_*\in\overline{\Omega}$ and $y_*\in\Omega$,
there exists a positive constant $C_3$ such that 
$$
u(x_*+\rho\nu,\rho^2)^p
\ge C_2^p v(x_*+\rho\nu,\rho^2)^p
\ge C_3\rho^{p(2\gamma q+2q+1)/q}
$$ 
for all sufficiently small $\rho>0$, 
where $\nu$ is as in Lemma~\ref{Lemma:5.1}. 
This together with \eqref{eq:5.14} implies that 
\begin{equation}
\label{eq:5.15}
\lim_{\rho\to 0+}\frac{\omega_{p,1/2}(\rho:x_*,y_*)}{\rho}=\infty. 
\end{equation}
Furthermore, 
it follows from \cite[Section~2]{Kenni} (see also property~(g) in Section~2) and \eqref{eq:5.14} that 
$$
g(x,t,v):=v^{3-1/p}t^{2\gamma}f(x)
$$
is $\beta$\,-concave in $D\times(0,\infty)$ with 
$$
\frac{1}{\beta}=3-\frac{1}{p}+2\gamma+\frac{1}{q}<1, 
$$
whence $g$ is concave in $D\times(0,\infty)$. 
Therefore, by \eqref{eq:5.12} and \eqref{eq:5.15} 
we apply Theorem~\ref{Theorem:3.1} with $\alpha=1/2$, and see 
that $u$ is parabolically $p$\,-concave in $\overline{D}$ in the case \eqref{eq:5.14}.  
Then assertion~(i) follows from property~(f) in Section~2 and Theorem~\ref{Theorem:4.1}. 
If $f$ is a positive constant function in $\Omega$, 
then $f$ is $q$\,-concave in $\Omega$ for any $q>1$. 
Therefore assertion~(ii) follows from assertion~(i) and property~(f) in Section~2, 
and the proof of Theorem~\ref{Theorem:5.1} is complete. 
$\Box$\vspace{5pt}
\newline
Next we state the following result on the optimality of 
the assumptions of Theorem~\ref{Theorem:5.1}. 
\begin{proposition}
\label{Proposition:5.1}
Let $\Omega$ be a bounded smooth convex domain in ${\bf R}^n$, $f$ positive smooth function in $\Omega$, 
and $1\le q\le\infty$. 
Assume that there exists a constant $C$ such that 
\begin{equation}
\label{eq:5.16}
f(x)\le 
\left\{
\begin{array}{ll}
C\,\mbox{{\rm dist}}(x,\partial\Omega)^{1/q} & \mbox{if}\quad q<\infty,\vspace{3pt}\\
C & \mbox{if}\quad q=\infty,
\end{array}
\right.
\end{equation}
for all $x\in\Omega$. 
Let $D:=\Omega\times(0,\infty)$ and $u$ satisfy \eqref{eq:5.1} with $0\le\gamma\le 1/2$. 
Then $u$ is not parabolically $p$\,-concave on $\overline{D}$ if 
\begin{equation}
\label{eq:5.17}
p>
\left\{
\begin{array}{ll}
\displaystyle{\frac{q}{1+2q+2\gamma q}} & \quad\mbox{for}\quad 1\le q<\infty,\vspace{5pt}\\
\displaystyle{\frac{1}{2(1+\gamma)}} & \quad\mbox{for}\quad q=\infty.
\end{array}
\right.
\end{equation}
\end{proposition}
{\bf Proof.}
Let $1\le q<\infty$, and assume \eqref{eq:5.17}.  
Let $\mu=\mu(\Omega)$ and $\psi>0$ be the first (positive) Dirichlet eigenvalue and eigenfunction (normalized, for instance, so that $\|\psi\|_{L^2(\Omega)}=1$) 
for $-\Delta$ on $\Omega$, respectively. 
Then there exists a positive constant $C_1$ such that 
\begin{equation}
\label{eq:5.18}
C_1^{-1}\mbox{{\rm dist}}(x,\partial\Omega)\le \psi(x)\le C_1\mbox{{\rm dist}}(x,\partial\Omega),\qquad x\in\Omega. 
\end{equation}
Put 
$$
w(x,t):=At^{\gamma+1}\psi(x)^{1/q},\qquad A>0.
$$
By \eqref{eq:5.16} and \eqref{eq:5.18} 
we can take a sufficiently large $A>0$ such that  
\begin{eqnarray*}
 & & \partial_t w-\Delta w-t^\gamma f\\
 & & =(\gamma+1)At^\gamma\psi(x)^{1/q}-\frac{A}{q}t^{\gamma+1}\psi^{1/q-1}\Delta\psi
-\frac{A}{q}\left(\frac{1}{q}-1\right)\psi^{1/q-2}|\nabla\psi|^2-t^\gamma f\\
 & & \ge\mu\frac{A}{q}t^{\gamma+1}\psi^{1/q}\ge 0\qquad\mbox{in}\quad D.
\end{eqnarray*}
Then the comparison principle together with \eqref{eq:5.18} yields
\begin{equation}
\label{eq:5.19}
0<u(x,t)\le w(x,t)\le C_2 t^{\gamma+1}\mbox{{\rm dist}}(x,\partial\Omega)^{1/q}
\quad\mbox{in}\quad D,
\end{equation}
for some constant $C_2>0$. 
Let $x_*\in\partial\Omega$ and $\nu_*$ be the inner unit normal vector to $\partial\Omega$ at $x_*$. 
We deduce from \eqref{eq:5.19} that 
$$
u(x_*+\rho\nu_*,\rho^2)^p\le C_2^p\rho^{p(2\gamma+2+1/q)}
$$
for all sufficiently small $\rho>0$. 
Since $p(2\gamma+2+1/q)>1$ by \eqref{eq:5.17}, 
$u(x_*+\rho\nu_*,\rho^2)^p$ is not concave with respect to $\rho$. 
This means that $u$ is not parabolically $p$\,-concave in $\overline{D}$ in the case $1\le q<\infty$. 
In the case $q=\infty$, by the comparison principle we see that 
$$
u(x,t)\le Bt^{\gamma+1}\quad\mbox{in}\quad D,
$$
for some positive constant $B$. 
Since 
$$
u(x,t^2)^p\le B^p t^{p(2\gamma+2)}\quad\mbox{in}\quad D, 
\qquad p(2\gamma+2)>1,
$$
$u$ is not parabolically $p$\,-concave in $\overline{D}$ in the case $q=\infty$. 
Thus Proposition~\ref{Proposition:5.1} follows.
$\Box$
\vspace{5pt}
\newline
Proposition~\ref{Proposition:5.1} means that 
in the case where $f$ is $q$\,-concave in $\Omega$ for some $q\ge 1$, 
under assumption~\eqref{eq:5.17}, 
the solution $u$ of \eqref{eq:5.1} is not necessarily parabolically $p$\,-concave in $\overline{D}$. 
Furthermore, we give additional comments on Theorem~\ref{Theorem:5.1} in the case $\gamma=0$. 
\begin{remark}
\label{Remark:5.1}
Let $u$ be a solution of problem~\eqref{eq:5.1} with $\gamma=0$.
\vspace{3pt}
\newline
{\rm (i)} It follows from Theorem~{\rm\ref{Theorem:5.1}} 
that if $f$ is $q$\,-concave in $\Omega$ for some $q\ge 1$, then 
the solution $u$ is parabolically $p$\,-concave in $\overline{D}$ 
with $p=q/(1+2q)$.  
Furthermore, the solution~$u$ converges to the unique solution $v$ of 
\begin{equation}
\label{eq:5.20}
\Delta v+f(x)=0\quad\mbox{in}\quad\Omega,
\qquad
v=0\quad\mbox{on}\quad\partial\Omega,
\end{equation}
pointwisely. 
On the other hand, 
if $u$ is parabolically $\tilde{p}$\,-concave in $\overline{D}$, 
then, due to property~{\rm(c)} in Section~{\rm 2}, 
the solution~$v$ of \eqref{eq:5.20} is $\tilde{p}$\,-concave in $\overline{\Omega}$. 
\vspace{3pt}
\newline
{\rm (ii)}
For any $q\in[1,\infty]$, 
there exists a $q$\,-concave function $f$ in $\Omega$ 
such that the solution $v$ of \eqref{eq:5.20} 
is not $r$\,-concave in $\overline{\Omega}$ for any $r>q/(1+2q)$
{\rm(}see Theorem~{\rm 6.2} in {\rm\cite{Kenni}}{\rm)}. 
In this case, assertion~{\rm (i)} implies that 
the solution $u$ is not parabolically $r$\,-concave in $\overline{D}$ 
for any $r>q/(1+2q)$.
\end{remark}

Next we deal with problem~\eqref{eq:1.1} 
in the case where $f(x,t,v,\theta)=v^\gamma$ with $\gamma\in(0,1)$. 
\begin{theorem}
\label{Theorem:5.2}
Let $\Omega$ is a bounded convex domain in ${\bf R}^n$ and $D:=\Omega\times(0,\infty)$. 
Consider the problem
\begin{equation}
\label{eq:5.21}
\partial_t u=\Delta u +u^\gamma\quad\mbox{in}\quad D,\qquad
u=0\quad\mbox{on}\quad\partial D,
\end{equation}
where $0<\gamma<1$. Then the maximal solution $u$ of \eqref{eq:5.21} is positive in $D$ 
and parabolically $p$\,-concave in $\overline{D}$ with 
\begin{equation}
\label{eq:5.22}
p=\frac{1-\gamma}{2}.
\end{equation}
Furthermore, $H(t)$ is $q$\,-concave in $(0,\infty)$ with 
$q=(1-\gamma)/[n(1-\gamma)+2]$.  
\end{theorem}
{\bf Proof.}
As in the proof of Theorem~\ref{Theorem:5.1}, 
it suffices to treat only the case where $\Omega$ is smooth. 
For any $\epsilon>0$, let $u_\epsilon$ be a solution of 
\begin{equation}
\label{eq:5.23}
\partial_t u=\Delta u +(u+\epsilon)^\gamma\quad\mbox{in}\quad D,\qquad
u=0\quad\mbox{on}\quad\partial D. 
\end{equation}
Then, by a similar argument to \eqref{eq:5.12} we have 
\begin{equation}
\label{eq:5.24}
\partial_t u_\epsilon\ge 0\quad\mbox{in}\quad D.
\end{equation} 
As before, let $\psi>0$ and $\mu$ be the first Dirichlet eigenfunction and eigenvalue for $-\Delta$ in $\Omega$, respectively. 
For $\delta>0$, set 
$$
w(x,t):=\delta t^{1/(1-\gamma)}\psi(x).
$$ 
Taking a sufficiently small $\delta>0$ if necessary, we have  
$$
\partial_t w-\Delta w-(w+\epsilon)^\gamma
\le \frac{1}{1-\gamma}\delta t^{\frac{\gamma}{1-\gamma}}\psi+\mu\delta t^{\frac{1}{1-\gamma}}\psi
-\delta^\gamma t^{\frac{\gamma}{1-\gamma}}\psi^\gamma\le 0
\quad\mbox{in}\quad\Omega\times(0,1)
$$
and 
$w(x,t)=0$ on $\partial\Omega\times(0,1)$ and $\overline{\Omega}\times\{0\}$. 
Then, applying the comparison principle, 
by \eqref{eq:5.18} we obtain 
\begin{equation}
\label{eq:5.25}
u_\epsilon(x,t)\ge w(x,t)\ge C_1\delta t^{\frac{1}{1-\gamma}}\mbox{{\rm dist}}(x,\partial\Omega)
\qquad\mbox{in}\quad\Omega\times[0,1),
\end{equation}
where $C_1$ is a positive constant independent of $\epsilon$. 

Let $\tilde{u}$ be a solution of \eqref{eq:5.21}. 
We apply the comparison principle again, and obtain  
$$
0\le \tilde{u}(x,t)\le u_{\epsilon_1}(x,t)\le u_{\epsilon_2}(x,t)\qquad\mbox{in}\quad\overline{D}
$$
if $0<\epsilon_1\le\epsilon_2$. 
This implies that the limit function 
$u(x,t):=\lim_{\epsilon\to 0}u_\epsilon(x,t)$
exists in $\overline{D}$. Furthermore, thanks to the regularity theorems for parabolic equations and \eqref{eq:5.24}, 
we see that $u$ is a solution of \eqref{eq:5.21} such that 
\begin{equation}
\label{eq:5.26}
0\le \tilde{u}(x,t)\le u(x,t)\quad\mbox{in}\quad\overline{D}
\qquad\mbox{and}\qquad
\partial_t u\ge 0\quad\mbox{in}\quad D. 
\end{equation}
In particular, we see that $u$ is the maximal solution of \eqref{eq:5.21}. 
Furthermore, by \eqref{eq:5.25} we obtain 
$$
u(x,t)\ge C_1\delta t^{\frac{1}{1-\gamma}}\mbox{{\rm dist}}(x,\partial\Omega)>0
\qquad\mbox{in}\quad\Omega\times(0,1),
$$
which implies that $u$ is a positive solution of \eqref{eq:5.21}. 

On the other hand, the function
\begin{equation}
\label{eq:5.27}
z(x,t):=u(x,t)^{1-\gamma} 
\end{equation}
satisfies 
\begin{equation}
\left\{
\begin{array}{l}
\partial_t z=\Delta z+\displaystyle{\left(\frac{1}{1-\gamma}-1\right)\frac{|\nabla z|^2}{z}}+1-\gamma
\ge\Delta z+1-\gamma\quad\mbox{in}\quad D,\vspace{7pt}\\
z=0\quad\mbox{on}\quad\partial D.
\end{array}
\right.
\label{eq:5.28}
\end{equation}
Let 
\begin{equation}
\label{eq:5.29}
\frac{1-\gamma}{3}<p<\frac{1-\gamma}{2}. 
\end{equation}
Let $U$ be the solution of \eqref{eq:5.2} with $\gamma=0$ and 
\begin{equation}
\label{eq:5.30}
0<d<\frac{1-\gamma}{p}-2.
\end{equation}
Then we apply the comparison principle to \eqref{eq:5.28} 
to obtain
$$
z(x,t)\ge C_2U(x,t)\quad\mbox{in}\quad D
$$
for some positive constant $C_2$. 
Furthermore, due to Lemma~\ref{Lemma:5.1}, 
for any $x_*\in\overline{\Omega}$ and $y_*\in\Omega$, 
there exists a positive constant $C_3$ such that 
$$
z(x_*+\nu\rho,\rho^2)\ge C_2U(x_*+\rho\nu,\rho^2)\ge C_3\rho^{d+2}
$$
for all sufficiently small $\rho>0$, where $\nu$ is as in Lemma~\ref{Lemma:5.1}. 
This together with \eqref{eq:5.27} implies 
$$
u(x_*+\rho\nu,\rho^2)^p\ge C_3^{p/(1-\gamma)}\rho^{p(d+2)/(1-\gamma)}
$$
for all sufficiently small $\rho>0$. 
Then we deduce from \eqref{eq:5.30} that 
\begin{equation}
\label{eq:5.31}
\lim_{\rho\to 0}\frac{\omega_{p,1/2}(\rho:x_*,y_*)}{\rho}=\infty. 
\end{equation}
Furthermore, by \eqref{eq:5.29} we see that 
the function 
$$
g(v):=v^{3-1/p}f(v^{1/p})=v^{3-(1-\gamma)/p}
$$
is concave with respect to $v\in(0,\infty)$. 
Therefore we apply Theorem~\ref{Theorem:3.1} with the aid of \eqref{eq:5.26} and \eqref{eq:5.31}, 
and see that for any $\lambda\in\Lambda_{n+1}$, 
$u_{1/2,p,\lambda}$ is a viscosity subsolution of \eqref{eq:5.21} such that $u_{1/2,p,\lambda}=0$ on $\partial D$. 
This implies that for any $\epsilon>0$,  $u_{1/2,p,\lambda}$ is a viscosity subsolution of \eqref{eq:5.23} 
such that $u_{1/2,p,\lambda}=0$ on $\partial D$. 
Since (WCP) holds for problem~\eqref{eq:5.23}, we have 
$$
u_{1/2,p,\lambda}(x,t)\le u_\epsilon(x,t)\qquad\mbox{in}\quad\overline{D}. 
$$
Therefore, letting $\epsilon\to 0$, we obtain  
$$
u_{1/2,p,\lambda}(x,t)\le u(x,t)\qquad\mbox{in}\quad\overline{D},
$$
and see that $u$ is parabolically $p$\,-concave in $\overline{D}$ with \eqref{eq:5.29}.
Hence, we deduce from properties~(e) and (f) in Section~2 
that $u$ is parabolically $p$\,-concave in $\overline{D}$ with $p=(1-\gamma)/2$. 
In addition, by Theorem~\ref{Theorem:4.1} we obtain the desired concavity property of $H(t)$, 
and the proof of Theorem~\ref{Theorem:5.2} is complete.
$\Box$
\vspace{3pt}
\newline
Similarly to Remark~\ref{Remark:5.1}~(i), 
the solution $u$ of \eqref{eq:5.21} converges to a positive solution $v$ of 
\begin{equation}
\label{eq:5.32}
\Delta v+v^\gamma=0\quad\mbox{in}\quad\Omega,
\qquad
v=0\quad\mbox{on}\quad\partial\Omega,
\end{equation}
pointwisely. Furthermore, it follows from Theorem~\ref{Theorem:5.2} 
that $v$ is $p$\,-concave in $\overline{\Omega}$ with $p=(1-\gamma)/2$. 
This coincides with the concavity property obtained by \cite[Theorem~4.2]{Kenni} 
for problem \eqref{eq:5.32} (see also \cite{LV2}). 
\bibliographystyle{amsplain}

\end{document}